\newif\ifarxived
\ifx\notArxived\Kaba\arxivedtrue\fi
\arxivedtrue
\newif\iftesting
\ifx\Testing\Kaba\else\testingtrue\fi
\newif\ifprivate
\ifx\Private\Kaba\else\privatetrue\fi
\ifarxived
\documentclass[12pt]{article}
\else
\documentclass[12pt,uplatex]{article}
\fi
\ifarxived
\usepackage[bookmarks=true,bookmarksnumbered=true,
  bookmarkstype=toc,
  pdftitle={The first-order definability of generic large cardinals},%
  pdfauthor={Sakaé Fuchino, and Hiroshi Sakai},%
  pdfsubject={},%
  pdfkeywords={generic supercompactness, Laver-generic supercompactness, saturated ideal, 
  elementary embedding}]{hyperref}
\else
\usepackage[dvipdfmx,bookmarks=true,bookmarksnumbered=true,
  bookmarkstype=toc,
  pdftitle={The first-order definability of generic large cardinals},%
  pdfauthor={Sakaé Fuchino, and Hiroshi Sakai},%
  pdfsubject={},%
  pdfkeywords={generic supercompactness, Laver-generic supercompactness, saturated ideal, 
  elementary embedding}]{hyperref}
\fi
\ifarxived
\else
\usepackage[dvipdfmx]{pxjahyper}
\hypersetup{
  colorlinks=true,
  linkcolor=blue,
  filecolor=blue,      
  citecolor=cyan,
  urlcolor=cyan
}
\fi
\ifarxived
\usepackage{graphicx, color}
\else
\usepackage[dvipdfmx]{graphicx, color}
\fi
\usepackage{amsmath, amssymb}
\usepackage{bbm}
\usepackage{dsfont}
\usepackage{bbold}
\usepackage{marginnote}
\usepackage{mathtools}
\usepackage[normalem]{ulem}
\usepackage{setspace}
\definecolor{darkelectricblue}{rgb}{0.33, 0.41, 0.47}
\newcommand{\extendedcolor}{\color{darkelectricblue}}
\newif\ifextended
\ifarxived\extendedtrue\else
\ifx\Extended\Kaba\else\extendedtrue\fi
\fi
\ifarxived\extendedtrue\fi
\definecolor{darkelectricblue}{rgb}{0.33, 0.41, 0.47} 
\newif\ifJapanese

\usepackage{theorem}
\newcommand{\bbd}[1]{{\mathbb{#1}}}
\theorembodyfont{\ifJapanese\rm\else\it\fi}
\newcount\minute	
\newcount\hour		
\newcount\hourMins  
\ifJapanese
\def\today%
{
  \the\year 年\,\zeroPadTwo{\the\month}月\,\zeroPadTwo{\the\day}日%
}
\fi
\def\now%
{
  \minute=\time    
  \hour=\time \divide \hour by 60 
  \hourMins=\hour \multiply\hourMins by 60
  \advance\minute by -\hourMins 
  \zeroPadTwo{\the\hour}:\zeroPadTwo{\the\minute}%
}
\def\zeroPadTwo#1%
{
  \ifnum #1<10 0\fi    
  #1
}
\date{}

\renewcommand{\baselinestretch}{1.2}
\renewcommand{\thefootnote}{(\arabic{footnote})\,}
\iftesting
\fi
\iftesting
\newcommand{\Label}[1]{\label{#1}\marginpar{{\renewcommand{\baselinestretch}{0.6}\tiny 
		  #1}}}
\else
\newcommand{\Label}[1]{\label{#1}}
\fi
\def\memo#1{\ifprivate\marginpar{{\normalsize\renewcommand{\baselinestretch}{0.4}\tiny%
			#1\par}}\else\fi}%

\newcounter{frml}[section]
\newcounter{frmla}[section]
\def\thefrml{{\arabic{section}.\arabic{frml}}}
\def\thefrmla{{$\aleph$\arabic{section}.\arabic{frmla}}}
\def\frmlabel#1{\refstepcounter{frml}{\def\baka{#1}\ifx\baka\empty\else\label{#1}\fi}%
{\rm({\thefrml})\hfill\hfill\hfill}}
\def\frmlabela#1{\refstepcounter{frmla}{\def\baka{#1}\ifx\baka\empty\else\label{#1}\fi}%
{\rm({\thefrmla})\hfill\hfill\hfill}}
\def\xitem[#1]{\item[\frmlabel{#1}]\mbox{}%
	\iftesting\marginpar{{\renewcommand{%
				\baselinestretch}{0.6}\tiny#1}}\fi\ignorespaces}
\def\xitemq[#1]{\item[\frmlabel{#1}]\mbox{}%
	\ignorespaces}
\def\xitema[#1]{\item[\frmlabela{#1}]\mbox{}%
	\iftesting\marginpar{{\renewcommand{%
				\baselinestretch}{0.6}\tiny#1}}\fi\ignorespaces}
\def\xitemsub[#1]#2{\item[\frmlabel{#1}$_{#2}$]\mbox{}%
	\iftesting\marginpar{{\renewcommand{%
				\baselinestretch}{0.6}\tiny#1}}\fi\ignorespaces}
\def\xitemcite[#1]{\item[\rlap{\rm(\ref{#1})}\hspace*{3em}]\mbox{}%
	\iftesting\marginpar{{\renewcommand{%
				\baselinestretch}{0.6}\tiny#1}}\fi\ignorespaces}
\def\xitemciteb[#1]#2{\item[\rlap{\rm(\ref{#1}#2)}\hspace*{3em}]\mbox{}%
	\iftesting\marginpar{{\renewcommand{%
				\baselinestretch}{0.6}\tiny#1}}\fi\ignorespaces}
\def\xxitem[#1][#2]{\item[(\ref{#1}{\makebox[1.4ex][c]{#2}})]\mbox{}%
	\iftesting\marginpar{{\renewcommand{%
				\baselinestretch}{0.6}\tiny\{#1\}\{#2\}}}\fi\ignorespaces}
\def\xitemof#1{{\rm({\ref{#1}})}}
\def\Xitem[#1]{\item[{\makebox[7ex][l]{\rm(\ref{#1})}}]\iftesting\marginnote{{\renewcommand{%
				\baselinestretch}{0.6}\tiny#1}}\fi\ignorespaces}

\newenvironment{xitemize}{\begin{list}{}{\parsep=0.5\smallskipamount%
			\itemindent=-0.4ex%
			\itemsep=0.5\smallskipamount\leftmargin=4em\labelwidth=3em\labelsep=0.7em}}%
							 {\end{list}}
\def\assert#1{\noindent\makebox[4.8ex][r]{\rm(\makebox[2.2ex][c]{#1})}\ \ \ignorespaces}
\def\assertof#1{\makebox[4ex][c]{\rm(\makebox[2.2ex][c]{#1})}}%
\def\wassert#1{\noindent\makebox[4.8ex][r]{\em{\rm(\makebox[2.3ex][c]{#1})}}\hspace{0.8em}\ignorespaces}
\def\wassertof#1{\makebox[4.2ex][c]{\rm(\makebox[2.3ex][c]{#1})}}%

\newcommand{\bysame}[1]{\underline{\phantom{#1}}}%

\newtheorem{Thm}{\ifJapanese{\bf 定理}\else {\bf Theorem}\fi}[section]

\newtheorem{ThmA}{\ifJapanese{\bf 定理\,A}\else{\bf Theorem\,A}\fi}[section]

\newtheorem{Prop}[Thm]{\ifJapanese{\bf 命題}\else{\bf Proposition}\fi}

\newtheorem{Lemma}[Thm]{\ifJapanese{\bf 補題}\else{\bf Lemma}\fi}
\newtheorem{LemmaA}[ThmA]{\ifJapanese{\bf 補題\,A}\else{\bf Lemma\,A}\fi}
\newtheorem{FactA}[ThmA]{Fact A}
\newtheorem{Cor}[Thm]{\ifJapanese{\bf 系}\else{\bf Corollary}\fi}

\newtheorem{Claim}{{\bf Claim}}[Thm]

\newtheorem{Subclaim}{{\bf Subclaim}}[Claim]

\newcommand{\prf}{\ifJapanese{\bf 証明．\ }\ignorespaces\else{\bf 
		Proof.\ \ }\ignorespaces\fi}

\newcommand{\prfofClaim}{\raisebox{-.4ex}{\Large $\vdash$\ \ }}

\newcommand{\prfof}[1]{\ifJapanese{\bf #1 の証明．\ \ }%
	\ignorespaces\else{\bf Proof of #1:}\ \ \ignorespaces\fi}
\newcommand{\Thmof}[1]{\ifJapanese{定理\,\ref{#1}}\else{Theorem~\ref{#1}}\fi}

\newcommand{\Lemmaof}[1]{\ifJapanese{補題\,\ref{#1}}\else{Lemma \ref{#1}}\fi}

\newcommand{\Propof}[1]{\ifJapanese{命題\,\ref{#1}}\else{Proposition~\ref{#1}}\fi}

\newcommand{\Claimof}[1]{{Claim \ref{#1}}}

\newcommand{\sectionof}[1]{\ifJapanese{第\ref{#1}節}\else{Section~\ref{#1}}\fi}

\newcommand{\footnoteof}[1]{\ifJapanese{脚注\ref{#1}}\else{footnote \ref{#1}}\fi}

\newcommand{\Claimabove}{{Claim \number\theClaim}}

\newsavebox{\qedbox}\sbox{\qedbox}{
{\unitlength=0.05mm \begin{picture}(40,60)
\put(0,0){\framebox(30,44)[cc]{}}
\put(30,-7){\rule{7\unitlength}{44\unitlength}}
\put(10,-7){\rule{27\unitlength}{7\unitlength}}
\end{picture}}}
\newcommand{\qed}{\mbox{}\hfill\usebox{\qedbox}}
\newcommand{\smallqed}%
{\mbox{}\smallskip\hfill\raisebox{-.4ex}{\Large $\dashv$}}
\newcommand{\qedof}[1]%
{\mbox{} \hspace*{\fill}{\usebox{\qedbox}{\tiny~(#1)}}}
\newcommand{\Qedof}[1]%
{\mbox{} \hspace*{\fill}{\usebox{\qedbox}%
{\tiny~(#1~\number\theThm)}}}
\newcommand{\QedAof}[1]%
{\mbox{} \hspace*{\fill}{\usebox{\qedbox}%
{\tiny~(#1~\number\theThmA)}}}
\newcommand{\qedofThm}{\Qedof{\ifJapanese 定理\else Theorem\fi}}

\newcommand{\qedofProp}{\Qedof{\ifJapanese 命題\else Proposition\fi}}
\newcommand{\qedofLemma}{\Qedof{\ifJapanese 補題\else Lemma\fi}}

\newcommand{\qedofFact}{\Qedof{Fact}}
\newcommand{\qedskip}{\medskip}

\newcommand{\qedofClaimmark}{\mbox{}\hfill\raisebox{-.4ex}{\Large $\dashv$}}
\newcommand{\qedofClaim}%
{\mbox{}\hfill\raisebox{-.4ex}{\Large $\dashv$ }\nolinebreak%
\mbox{\tiny~(Claim~\number\theClaim)}}
\newcommand{\qedofSubclaim}%
{\mbox{}\hfill\raisebox{-.4ex}{\Large $\dashv$ }\nolinebreak%
\mbox{\tiny~(Subclaim~\number\theSubclaim)}}
\ifarxived
\newcommand{\ubecause}[3]{\underbrace{{}#1{}%
  \ifx\bakakaba#2\bakakaba\rule[-0.72ex]{0pt}{1pt}\else\rule[#2]{0pt}{1pt}\fi}_{\mbox{\footnotesize\clap{#3}}}}
\newcommand{\obecause}[3]{\overbrace{{}#1{}%
  \ifx\bakakaba#2\bakakaba\rule[1.62ex]{0pt}{1pt}\else\rule[#2]{0pt}{1pt}\fi}^{\mbox{\footnotesize\clap{#3}}}}
\else
\newcommand{\ubecause}[3]{\underbrace{{}#1{}%
  \ifx\bakakaba#2\bakakaba\rule[-0.4zh]{0pt}{1pt}\else\rule[#2]{0pt}{1pt}\fi}_{\mbox{\footnotesize\clap{#3}}}}
\newcommand{\obecause}[3]{\overbrace{{}#1{}%
  \ifx\bakakaba#2\bakakaba\rule[0.9zh]{0pt}{1pt}\else\rule[#2]{0pt}{1pt}\fi}^{\mbox{\footnotesize\clap{#3}}}}
\fi
\newcommand{\cardof}[1]{\mathopen{|\,}#1\mathclose{\,|}}

\newcommand{\setof}[2]{\{#1\,:\,#2\}}
\newcommand{\ssetof}[1]{\{#1\}}

\newcommand{\subseteqand}[1]{\mathrel{\mathop{\subseteq}%
		\limits_{\scriptscriptstyle\hbox to 14pt{$\scriptscriptstyle #1$\hss}}}}

\newcommand{\mapping}[3]{#1:#2\rightarrow #3}

\newcommand{\elembed}[3]{#1:#2\stackrel{\preccurlyeq\hspace{0.8ex}}{\rightarrow}#3}
\newcommand{\combed}[3]{#1:#2\stackrel{{\scriptstyle\mathrel{{\leqslant}%
		\hspace{-0.72ex}{\lower-0.27ex\hbox{\scalebox{0.7}{$\scriptstyle\circ$}}}}\hspace{0.7ex}}}{\rightarrow}#3}

\newcommand{\imageof}{{}^{\,{\prime}{\prime}}}
\newcommand{\seqof}[2]{\langle#1\,:\,#2\rangle}
\newcommand{\pairof}[1]{\langle#1\rangle}

\newcommand{\forces}[2]{\,\|\hspace{-.35ex}\mbox{\sf--}_{\,#1\,}%
\mbox{\rm``}\,#2\,\mbox{\rm''}}
\newcommand{\forcessanf}[2]{\,\|\hspace{-.35ex}\mbox{\sf--}_{\,#1\,}%
\mbox{\rm``}\,#2}

\newcommand{\modelof}[1]{\models\!\mbox{\rm``\,}#1\mbox{\rm''}}

\newcommand{\crit}{\mbox{\it crit\/}}

\newcommand{\circleq}{\mathrel{{\leqslant}%
		\hspace{-0.86ex}{\lower-0.53ex\hbox{$\scriptscriptstyle\circ$}}}}

\newcommand{\variables}[2]{{#1}_0\ctentenc {#1}_{#2}}

\newcommand{\restr}{\restriction}

\newcommand{\Col}{{\rm Col}}

\newcommand{\otp}{\mathop{\mbox{\it otp\/}}}
\newcommand{\dom}{\mathop{\rm dom}}

\newcommand{\id}{{\rm id}}

\newcommand{\trcl}{{\it trcl\/}}

\newcommand{\poP}{\bbd{P}}

\newcommand{\genG}{\mathbb{G}}
\newcommand{\genH}{\mathbb{H}}

\newcommand{\condp}{\mathbbm{p}}
\newcommand{\condq}{\mathbbm{q}}

\newcommand{\LT}{{<}\,}

\newcommand{\ctenten}{,\mbox{}\hspace{0.08ex}{.}{.}{.}\hspace{0.1ex}}
\newcommand{\ctentenc}{,{}\linebreak[0]\hspace{0.04ex}{{.}{.}{.}\hspace{0.1ex},\,}\linebreak[0]}

\newcommand{\threedotsh}{\hspace{0.14ex}{\cdot}{\cdot}{\cdot}\hspace{0.1ex}}

\newcommand{\calB}{{\mathcal B}}

\newcommand{\calF}{{\mathcal F}}

\newcommand{\calH}{{\mathcal H}}
\newcommand{\calI}{{\mathcal I}}
\newcommand{\calL}{{\mathcal L}}

\newcommand{\calP}{{\mathcal P}}
\newcommand{\calQ}{{\mathcal Q}}

\newcommand{\calW}{{\mathcal W}}

\newcommand{\utf}{\utilde{f}}

\newcommand{\utpoQ}{\utilde{\mathbb Q}}

\newcommand{\Lin}{{\calL}_{\varin}}
\newcommand{\varin}{\mathrel{\,\varepsilon\,}}

\newcommand{\ZF}{{\sf ZF}}
\newcommand{\ZFC}{{\sf ZFC}}

\newcommand{\NBGC}{{\sf NBGC}}

\newcommand{\PFA}{{\sf PFA}}

\newcommand{\st}{such that}
\newcommand{\wrt}{with respect to}

\newcommand{\tfae}{the following are equivalent}

\newcommand{\uniV}{{\sf V}}
\newcommand{\uniW}{{\sf W}}
\newcommand{\po}{poset}

\newcommand{\pos}{posets}

\newcommand{\Ba}{Boolean algebra}
\newcommand{\Pkl}[2]{\ifx\bakakaba#1\bakakaba\ifx\bakakaba#2\bakakaba{\mathcal 
    P}_\kappa(\lambda)\else{\mathcal P}_\kappa(#2)\fi\else{\mathcal P}_{#1}(#2)\fi}

\newcommand{\utildeT}[1]{%
	\hbox to 0pt{$\mathop{#1}\limits_{\raise0.25ex\hbox{$\scriptstyle\sim$}}$\hss}%
		\relax\phantom{\underline{#1}}}
\newcommand{\utildeS}[1]{%
	\hbox to 0pt{\smash{$\mathop{\scriptstyle #1}\limits_{%
				\raisebox{0.6ex}[0pt]{$\scriptscriptstyle\sim$}}$}\hss}%
		\relax\phantom{\mathord{{#1}_{\rule[-0.6ex]{0pt}{1pt}}}}}
\newcommand{\utildeSS}[1]{%
	\hbox to 0pt{$\mathop{\scriptscriptstyle #1}%
		\limits_{\scriptscriptstyle\sim}$\hss}%
		\relax\phantom{\underline{#1}}}
\newcommand{\utilde}[1]{%
	\mathchoice{\utildeT{#1}}{\utildeT{#1}}{\utildeS{#1}}{\utildeSS{#1}}}

{\end{minipage}\end{trivlist}}

\begin{document}
\addcontentsline{toc}{section}{* The first-order definability of generic large cardinals}
\title{\bf The first-order definability of generic large cardinals}
\ifarxived
\author{\bf Saka\'e Fuchino (\quad\quad \quad), and Hiroshi Sakai (\quad\quad \quad\quad)}
\else
\author{\bf Saka\'e Fuchino (渕野 昌), and Hiroshi Sakai (酒井 拓史)}
\fi
\maketitle
\renewcommand{\thefootnote}{$\ast$\ }
  \footnotetext{Graduate School of System Informatics, Kobe University \\Rokko-dai 1-1, Nada, Kobe 657-8501 Japan
   \\
    \scalebox{0.95}[1]{\tt fuchino\@diamond.kobe-u.ac.jp,\,hsakai@people.kobe-u.ac.jp}}


\begin{abstract}
We show that the notions of generic and Laver-generic supercompactness are first-order 
definable in the language of \ZFC. This also holds for generic and Laver-generic (almost) hugeness as 
well as for generic versions of other large cardinals. 
\end{abstract}


\ifextended
{\color{darkelectricblue}

\begin{quote}
	\footnotesize
	\noindent
	\centerline{\normalsize\tt Contents\hspace{6em}\mbox{}}\mbox{}\\
       {\mbox{}\hspace{-1.6em}\tt\makebox[3.4ex][l]{\ref{defbility:intro}.}%
         Introduction}\ \ \dotfill\ \ {\pageref{defbility:intro}}\\ 
       {\mbox{}\hspace{-1.6em}\tt\makebox[3.4ex][l]{\ref{defbility:normal}.}%
         $\uniV$-normal ultrafilters}\ \ \dotfill\ \ {\pageref{defbility:normal}}\\ 
       {\mbox{}\hspace{-1.6em}\tt\makebox[3.4ex][l]{\ref{defbility:initial}.}%
         Sufficiently large initial segment of elementary embeddings}\ \ \dotfill\ \ {\pageref{defbility:initial}}\\ 
       {\mbox{}\hspace{-1.6em}\tt
         References}\ \ \dotfill\ \ {\pageref{defbility:ref}}\\ 
\end{quote}}
\fi
\renewcommand{\thefootnote}{}
\footnotetext{{\it Date:} June 27, 2021
  \qquad {\it Last update:} 
  \today\ (\now\ \ifarxived WGST\else JST\fi)\vspace{-1\smallskipamount}
}
\footnotetext{\vspace{-1\smallskipamount}}
\footnotetext{The first author was partially supported by JSPS Kakenhi Grant No.\ 20K03717.
  The second author is 
  supported by JSPS Kakenhi Grant No.\ 18K03397.
\ifextended\\ \par
  {\color{darkelectricblue} This is an extended version of the paper with the same title. Some extra remarks and 
    details omitted in the final version for publication, as well as further corrections 
    after the publication, may be found 
    in this version. The additional 
    stuff is typeset in dark electric blue like this paragraph. The most recent edition of 
    this extended version is downloadable as:\smallskip\\ \qquad 
    \url{https://fuchino.ddo.jp/papers/definability-of-glc-x.pdf}}\else An extended version of the paper 
  with some details omitted in the submitted version, and with some possible further update is 
  downloadable as: \url{https://fuchino.ddo.jp/papers/definability-of-glc-x.pdf} 
\fi}

\renewcommand{\thefootnote}{(\arabic{footnote})\,}
\section{Introduction}\Label{defbility:intro}
For a class of \pos\ $\calP$, 
a cardinal $\kappa$ is said to be {\it generically supercompact by $\calP$} if,  
for any regular $\lambda\geq\kappa$, there is a \po\ $\poP\in\calP$ \st, for 
a $(\uniV,\poP)$-generic $\genG$, there are $M$, $j\subseteq\uniV[\genG]$ \st\ 
\begin{xitemize}
\xitem[def:x-intro-a] 
  $\elembed{j}{\uniV}{M\subseteq\uniV[\genG]}$, \footnote{When we write
    $\elembed{j}{\uniV}{M\subseteq\uniV[\genG]}$, we assume that $M$ is a transitive class in
    (and thus an inner model of) $\uniV[\genG]$.}
\xitem[def:x-intro-0] $\crit(j)=\kappa$, $j(\kappa)>\lambda$,
\xitem[def:x-intro-1] $j\imageof\lambda\in M$. 
\end{xitemize}
We shall call the class mapping $j$ as above a $\lambda$-generically supercompact embedding 
for $\kappa$ (in $\uniV[\genG]$). 

It is easy to see that a generically supercompact cardinal $\kappa$ for any class $\calP$ 
of \pos\ is regular. \memo{talk-nagoya-2019-05-31.pdf}
Even so, a generically supercompact cardinal can be a successor cardinal: If we 
collapse all cardinals below a supercompact cardinal $\kappa$ 
by $\Col(\omega_1,\kappa)$,\footnote{We use here Kanamori's notation in \cite{defbility:bib-kanamori} of L\'evy collapse.}
in the generic extension, $\kappa=\aleph_2$ and $\aleph_2$ is generically supercompact by
$\sigma$-closed \pos. 

A generically supercompact cardinal can $\kappa$ be also weakly inaccessible. Actually 
$\kappa$ can be even really supercompact for any $\calP$ as far as this $\calP$ contains 
the trivial \po. However, a generically supercompact $\kappa$ can also be weakly 
inaccessible (and much more) while it is not strongly inaccessible: If $\kappa$ is 
supercompact and $\kappa$ many Cohen reals are added, then  
$\kappa$ is still a regular inaccessible cardinal (and actually much more) and it remain 
generically supercompact by c.c.c.\ \pos, while it is 
the continuum in the generic extension.  

Similarly to the genuine supercompactness, it is not immediately clear if the notion of 
generic supercompactness is definable in 
the language of \ZFC. In most of the cases, this does not bother. This is because the generically 
supercompactness may be used in many applications merely as a schematic framework in which arguments in 
different settings are put together to obtain a better perspective. 

However, the circumstances become different if we would like to think generic supercompactness as a set-theoretic 
axiom. 

In \cite{defbility:bib-koenig}, Bernhard K\"onig gave a characterization of the statement
``$\omega_2$ is generically supercompact by $\sigma$-closed \pos'' in terms of the reflection 
of the non-existence of winning strategy of the second player in certain type of two player games. 
Since this reflection principle which K\"onig called ``Strong Game Reflection Principle'' 
is first-order definable, the statement mentioned above is also first-order formalizable.

In \cite{defbility:bib-I}, K\"onig's characterization is generalized to a characterization 
of the statement ``$\kappa^+$ is generically supercompact for $\LT\kappa$-closed forcing'' 
for arbitrary regular uncountable $\kappa$. By the same argument as above, we conclude 
from this result that the statement is also first-order formalizable. 

Based on the main idea in the proof of these results,  
we show 
in the following \sectionof{defbility:normal} that the generically supercompactness 
for any class $\calP$ of \pos\ is first-order definable. 

We say that a class $\calP$ of \pos\ {\it iterable}, if $\calP$ is closed \wrt\ 
restriction (i.e., if $\poP\in\calP$ and $\condp\in\poP$, 
then $\poP\restr\condp\in\calP$)\footnote{For the use of this condition, see the argument around \xitemof{def:x-gen-sc-0}}, and, 
for any $\poP\in\calP$ and $\poP$-name $\utpoQ$, we have 
\begin{xitemize}
\item[]  if $\forces{\poP}{\utpoQ\varin\calP}$ then $\poP\ast\utpoQ\in\calP$. 
\end{xitemize}

For a cardinal $\kappa$ and an iterable class $\calP$ of \pos, we call $\kappa$ a 
{\it Laver-generically supercompact for $\calP$} (or {\it L-g supercompact}, for short) if, 
for any $\lambda\geq\kappa$ and any 
$\poP\in\calP$, there is a $\poP$-name of a \po\ $\utpoQ$ with $\forces{\poP}{\utpoQ\varin\calP}$
\vspace{\smallskipamount}\st, for any $(\uniV,\poP\ast\utpoQ)$-generic filter $\genH$,  there 
are $M$, $j\subseteq\uniV[\genH]$ \st\ 
\begin{xitemize}
\xitem[defbility:x-gen-sc-13] $\elembed{j}{\uniV}{M}$, 
\xitem[defbility:x-gen-sc-14] $\crit(j)=\kappa$, $j(\kappa)
  >\lambda$,
\xitem[defbility:x-gen-sc-15] $\poP$, $\genH\in M$ and 
\xitem[defbility:x-gen-sc-16] $j\imageof\lambda\in M$.
\end{xitemize}

We shall call $j$ as above {\it a $\lambda$ L-g supercompact embedding (with the critical 
  point $\kappa$, associated with $\genH$ over $\uniV$)}. 

For $\calP=$ all the $\sigma$-closed \pos, the supercompact $\kappa$ in the ground model 
collapsed to be $\aleph_2$ by $\Col(\omega_1,\kappa)$ is L-g supercompact for $\calP$. 
For $\calP=$ all the proper \pos, the continuum in the standard model of \PFA\ obtained by 
starting from a supercompact $\kappa$ and by iterating with proper \pos\ with countable 
support along with a Laver diamond is L-g supercompact for $\calP$.

In these two models the L-g supercompact cardinal is $\aleph_2$. This is not a coincidence: 
If all elements of $\calP$ preserves $\omega_1$ and
$\Col(\omega_1,\ssetof{\omega_1})\in\calP$ then $\kappa$ being L-g supercompact for $\calP$ 
implies $\kappa=\aleph_2$ (\cite{defbility:bib-II}).

For $\calP=$ all the ccc \pos, a L-g supercompact cardinal for $\poP$ is obtained by 
starting from a supercompact $\kappa$ and then iterating $\kappa$-times by ccc \pos\ with finite support 
along with a Laver diamond.

The method in \sectionof{defbility:normal} cannot be applied (at least not in a straightforward 
way) to show the definability of 
Laver-generic large cardinals since apparently it cannot cover the condition 
\xitemof{defbility:x-gen-sc-15}.

In \sectionof{defbility:initial}, we show that the 
existence of generic elementary embedding can be recovered from a large enough initial 
segment of a generic elementary embedding (\Propof{def:P-is-1}). Using this, we can establish the 
definability of Laver-generic supercompactness for any iterable class of \pos (\Thmof{def:P-is-2}). 

The results discussed in this paper can be easily modified to adopt to other generic and 
Laver-generic large cardinals like those corresponding to super almost huge or super-huge 
cardinals. 

In the following, we assume that our formal framework is that of \ZFC\ and $\Lin$ denotes 
the language of set theory with the sole binary relation symbol $\varin$. Nevertheless, when we 
consider generic elementary embeddings which may not be first-order definable, we go over 
to the second-order framework of the axiom system of von Neumann-Bernays-Gödel (\NBGC) 
e.g.\ by adding an appropriate axiom $\Psi$ claiming the existence of certain (class) names 
of elementary embeddings in a generic extension over each \pos\ in a given class of \pos. 

We say that such system is first-order definable if we can find an axiom $\psi$
in $\Lin$ \st\ the original second-order axiom \NBGC\ $+$ $\Psi$ is a conservative extension of the the axiom 
system $\ZFC$ $+$ $\psi$. 

In the framework of \ZFC, when we are talking about a class $\calP$ of \pos, we assume that 
we fix an $\Lin$-formula $P(\cdot)$ which describes the elements of $\calP$ in such a way 
that $\calP=\setof{\poP}{P(\poP)}$. In this respect, when we 
said $\forces{\poP}{\utpoQ\varin\calP}$ in connection with iterability of $\calP$ above, we actually 
meant $\forces{\poP}{P(\utpoQ)}$. 

\section{$\uniV$-normal ultrafilters}
\Label{defbility:normal}
\ifextended{\extendedcolor
In the context of generic supercompactness, the condition \xitemof{def:x-intro-1}
implies a certain kind of closedness of $M$. This can be seen in 
the following Lemma:

\begin{LemmaA}{\rm (Lemma 2.5 in \cite{defbility:bib-II})}
  \Label{def:P-gen-sc-0}
  Suppose that $\genG$ is a $(\uniV,\poP)$-generic filter for a
  \po\/ $\poP\in\uniV$, and $\elembed{j}{\uniV}{M\subseteq\uniV[\genG]}$ is \st, for 
  cardinals $\kappa$, $\lambda$ in $\uniV$ 
  with $\kappa\leq\lambda$, $\crit(j)=\kappa$ and $j\imageof\lambda\in M$. Then, we have 
  the following:\smallskip

  \assert{1} For any set $A\in\uniV$  
  with $\uniV\models\cardof{A}\leq\lambda$, we have $j\imageof A\in M$. \smallskip

  \assert{2} $j\restr\lambda$, $j\restr\lambda^2\in M$.\smallskip

  \assert{3} For any $A\in\uniV$ with $A\subseteq\lambda$ or $A\subseteq\lambda^2$ 
  we have $A\in M$.\smallskip

  \assert{4} $(\lambda^+)^M\geq(\lambda^+)^\uniV$, Thus, if
  $(\lambda^+)^\uniV=(\lambda^+)^{\uniV[\genG]}$,  then
  $(\lambda^+)^M=(\lambda^+)^\uniV$. \smallskip

  \assert{5} $\calH(\lambda^+)^\uniV\subseteq M$.\smallskip

  \assert{6} $j\restr A\in M$ for all $A\in\calH(\lambda^+)^\uniV$. 
  \qed
\end{LemmaA}

In the following, we use Kanamori's notation of collapsing \pos\  
(see \S10 of \cite{defbility:bib-kanamori}). 

As it is already noticed in the introduction, 
it is consistent (modulo a supercompact cardinal) that a successor cardinal of a regular 
uncountable cardinal is generically 
supercompact. 

\begin{FactA}
  \Label{def:P-gen-sc-0-0}
  Suppose that $\kappa$ is a (really) supercompact cardinal, $\mu<\kappa$ a 
  regular uncountable cardinal, and $\poP_0=\Col(\mu,\kappa)$.
  Then, for 
  a $(\uniV,\poP_0)$-generic $\genG_0$, 
  \begin{xitemize}
  \item[] 
    $\uniV[\genG_0]\modelof{\mu^+\mbox{ is a generically supercompact cardinal by }
      \LT\mu\mbox{-closed \pos\,}}$.
  \end{xitemize}
\end{FactA}
\prf Note that $\uniV[\genG_0]\modelof{\mu^+=\kappa}$. 

For $\lambda\geq\kappa$, let $\elembed{j}{\uniV}{M}$ be a $\lambda$-supercompact 
embedding for $\kappa$. Then we have
\begin{xitemize}
\item[] $j(\poP_0)\ubecause{=}{}{by elementarity\qquad}\Col(\ubecause{j(\mu),}{}{\qquad$=\mu
  $}j(\kappa))^M\obecause{=}{}{\qquad by closedness of $M$}\Col(\mu,j(\kappa))^\uniV$. 
\end{xitemize}

For a $(\uniV[\genG_0], \Col(\mu,j(\kappa)\setminus\kappa))$-generic filter $\genG$,  
the lifting 
\begin{xitemize}
\item[] 
  $\elembed{\tilde{j}}{\uniV[\genG_0]}{\ubecause{M[\genG_0][\genG]}{}{$\subseteq\uniV[\genG_0][\genG]$}}$;\quad
  $\utilde{a}^{\genG_0}\mapsto j(\utilde{a})^{\genG_0\ast\genG}$
\end{xitemize}
witnesses the 
generic 
$\lambda$-supercompactness of $\ubecause{\kappa}{}{\qquad$=(\mu^+)^{\uniV[\genG_0]}$}$ by $\mu$-closed 
\pos\ in $\uniV[\genG_0]$.\ifarxived\vspace{-1.8ex}\else\vspace{-1zh}\fi\qedofFact\qedskip
}\fi 

For a class $\calP$ of \pos\ \st\ no $\poP\in\calP$ adds any 
new $\omega$-sequence of ground model sets, the first-order definability of the generic 
supercompactness by $\calP$ can be seen in the following Proposition. The 
Proposition can be shown by a direct imitation of the proof of the characterization of 
supercompactness by Solovay and Reinhardt in terms of the existence of normal ultrafilters 
(see e.g.\ Theorem 22.7 in \cite{defbility:bib-kanamori}).  

\begin{Thm}
\Label{def:P-gen-sc-1} Suppose that $\calP$ is a class of \pos\ \st\ no $\poP\in\calP$ adds any 
new $\omega$-sequence of ground model sets, and $\calP$ is closed \wrt\ restriction (i.e, 
if $\poP\in\calP$ and $\condp\in\poP$, then $\poP\restr\condp\in\calP$). 

An uncountable cardinal $\kappa$ 
is generically supercompact by $\calP$ if and only if,  
for any $\lambda\geq\kappa$, there is a $\poP\in\calP$ \st\ 
\begin{xitemize}
\item[] 
  $\forces{\poP}{\mbox{there is a }\uniV\mbox{-normal ultrafilter on }
  \calP^\uniV(\Pkl{}{}^\uniV)}.$
\end{xitemize}
\end{Thm}

Here, the notion of $\uniV$-normal ultrafilter is defined as follows:\quad 
Suppose that we are living in a universe $\uniW$ and $\uniV$ is an inner model in $\uniW$. 
Let $\lambda$ be an ordinal in $\uniV$, $\calI\in\uniV$, 
$\calI\subseteq\calP^\uniV(\lambda)$ a $\sigma$-ideal with $\ssetof{\xi}\in\calI$ for all
$\xi<\lambda$, and $\calB\in\uniV$ the  sub-\Ba\ 
$\calB=\calP^\uniV(\calI)$ of $\calP^\uniW(\calI)$. 

In $\uniW$, $U\subseteq\calB$ is a {\it$\uniV$-normal ultrafilter}
if 
\begin{xitemize}
\xitem[defbility:x-gen-sc-3] $U$ is a ultrafilter on the \Ba\ $\calB$. I.e., 
  \begin{xitemize}
  \item[\wassertof{i}] $\emptyset\not\in U$;
  \item[\wassertof{ii}] $A\cap A'\in U$ for any $A$, $A'\in U$;
  \item[\wassertof{iii}] if $A\in U$, $A\subseteq A'\in\calB$, then $A'\in U$; and
  \item[\wassertof{iv}] for any
    $A\in\calB$, either $A\in U$ or $\calI\setminus A\in U$; 
  \end{xitemize}
\xitem[defbility:x-gen-sc-4] 
  For any $x_0\in\calI$, we have 
  $\setof{x\in\calI}{x_0\subseteq x}\in U$;
\xitem[defbility:x-gen-sc-5] 
  For any $\seqof{A_\xi}{\xi\in\lambda}\in\uniV$, if
  $\setof{A_\xi}{\xi<\lambda}\subseteq U$, we have\\
  $\triangle_{\xi\in\lambda}A_\xi\in U$. Here, $\triangle_{\xi\in\lambda}A_\xi$ is the 
  diagonal intersection of $A_\xi$'s defined by 
  \begin{xitemize}
  \xitem[defbility:x-gen-sc-5-a] 
    $\triangle_{\xi\in\lambda}A_\xi:=\setof{x\in\calI}{x\in A_\xi
      \mbox{ for all }\xi\in x}$.
  \end{xitemize}
\end{xitemize}

\begin{Lemma}
  \Label{def:P-gen-sc-2} Suppose that $U\subseteq\calB$ is a
  $\uniV$-normal ultrafilter. \smallskip

  \wassert{1} For $\delta<\lambda$ \st\
  $\delta\in\calI$, and $\seqof{A_\xi}{\xi\in\delta}\in\uniV$ with 
  $A_\xi\in U$ for all $\xi\in\delta$, we have $\bigcap_{\xi\in\delta}A_\xi\in U$. 
  \smallskip

  \wassert{2} {\rm(Pressing Down Lemma)} For any $f\in\uniV$ with $\mapping{f}{\calI}{\uniV}$, if\/
  $\setof{x\in\calI}{f(x)\in x}\in U$, then there is $\xi<\lambda$ \st\
  $\setof{x\in\calI}{f(x)=\xi}\in U$.
\end{Lemma}
\prf  \assertof{1}:\quad Let $A_\xi:=\calI$ for all $\xi\in\lambda\setminus\delta$. Then
\begin{xitemize}
\item[] 
  $\ubecause{
    \obecause{\triangle_{\xi\in\lambda}A_\xi}{}{\qquad\qquad$\in U$ by \xitemof{defbility:x-gen-sc-5}}
    \cap
    \ubecause{\setof{x\in\calI}{\delta\subseteq x}}{}{\qquad$\in U$ by \xitemof{defbility:x-gen-sc-4}}
    }{}{\qquad\qquad$\in U$ by \xitemof{defbility:x-gen-sc-3},\,\assertof{ii}}
    \subseteq\bigcap_{\xi\in\delta}A_\xi$.
\end{xitemize}
Hence, $\bigcap_{\xi\in\delta}A_\xi\in U$ by \xitemof{defbility:x-gen-sc-3},\,\assertof{iii}.\smallskip\smallskip

\assertof{2}:\quad Suppose that $f$ is a counter-example to the assertion. That is,
\begin{xitemize}
\xitem[defbility:x-gen-sc-5-0] $A:=\setof{x\in\calI}{f(x)\in x}\in U$, but
\xitem[defbility:x-gen-sc-5-1] $A_\xi:=\setof{x\in\calI}{f(x)\not=\xi}\in U$ for all $\xi\in\lambda$.
\end{xitemize}
Then $\triangle_{\xi<\lambda}A_\xi\cap A\in U$ by \xitemof{defbility:x-gen-sc-5} and 
\xitemof{defbility:x-gen-sc-3},\,\assertof{ii}. By \xitemof{defbility:x-gen-sc-3},\,\assertof{i}, there is an 
element $x^*$ of this set. $f(x^*)\in x^*$ by \xitemof{defbility:x-gen-sc-5-0} but $f(x^*)\not=\xi$ for all 
$\xi\in x^*$ by \xitemof{defbility:x-gen-sc-5-1} and the definition \xitemof{defbility:x-gen-sc-5-a} of 
diagonal intersection. This is a contradiction.
\qedofLemma
\qedskip

\noindent
\prfof{\Thmof{def:P-gen-sc-1}} 
``$\Rightarrow$'':\quad 
Let $\lambda\geq\kappa$ and let $\poP$ be a $\LT\mu$-closed \po\ 
with $(\uniV,\poP)$-generic  $\genG$ and classes $j$, $M\subseteq\uniV[\genG]$ \st\
$\elembed{j}{\uniV}{M}$ is a $\lambda$-generically supercompact embedding for $\kappa$.
In particular, we have $j\imageof{\lambda}\in M$. Note that 
\begin{xitemize}
\xitem[defbility:x-gen-sc-5-2] 
  $M\models j\imageof{\lambda}\in\Pkl{j(\kappa)}{j(\lambda)}=j(\Pkl{}{}^\uniV)$. 
\end{xitemize}

In $\uniV[\genG]$, let 
\begin{xitemize}
\xitem[defbility:x-gen-sc-5-2-a] $U_j:=\setof{A\in \uniV}{A\subseteq\Pkl{}{}^\uniV,\, j\imageof{\lambda}\in j(A)}$. 
\end{xitemize}

\ifextended\else The following is easy to check:\fi
\begin{Claim}
  \Label{Cl-gen-sc-0}
  $U_j$ is a $\uniV$-normal ultrafilter on $\calP^\uniV(\Pkl{}{}^\uniV)$. 
\end{Claim}
\prfofClaim 
$U_j\models\xitemof{defbility:x-gen-sc-3}$, 
\assertof{i}:\quad  $j(\emptyset)=\emptyset$ by elementarity (and transitivity of $M$). Thus
$\emptyset\not\in U_j$ by definition. 

\assertof{ii}:\quad Suppose $A$, $A'\in U_j$. By definition 
this means that $j\imageof{\lambda}\in j(A)$ and $j\imageof{\lambda}\in j(A')$. It 
follows that
$j\imageof{\lambda}\in j(A)\cap j(A')\ubecause{=}{}{by elementarity}j(A\cap A')$. This 
shows that $A\cap A'\in U_j$. 

\assertof{iii}: Suppose that $A\in U_j$ and $A'\in\uniV$ is \st\
$A\subseteq A'\subseteq\Pkl{}{}^\uniV$. Then by elementarity we have
$M\models j(A)\subseteq j(A')$. Hence $j\imageof{\lambda}\in j(A)\subseteq j(A')$, and
$A'\in U_j$. 

\assertof{iv}: If $A\in\calP^\uniV(\Pkl{}{}^\uniV)\setminus U_j$, then by 
\xitemof{defbility:x-gen-sc-5-2},
$j\imageof{\lambda}\in j(\Pkl{}{}^\uniV)\setminus j(A)=j(\Pkl{}{}^\uniV\setminus A)$.
Thus $\Pkl{}{}^\uniV\setminus A\in U_j$. \smallskip

$U_j\models\xitemof{defbility:x-gen-sc-4}$:  Suppose $x_0\in\Pkl{}{}^\uniV$ and let
$A:=\setof{x\in\Pkl{}{}^\uniV}{x_0\subseteq x}$. 
Clearly $A\in\calP^\uniV(\Pkl{}{}^\uniV)$. By elementarity, and noting that
$j(x_0)=j\imageof{x_0}$ since $\cardof{x_0}<\kappa$, we have
\begin{xitemize}
\item[] $M\models j(A)=\setof{x\in\Pkl{j(\kappa)}{j(\lambda)}}{\ubecause{j(x_0)}{}{
$=j\imageof{x_0}$}\!\subseteq x}$.
\end{xitemize}
Thus $M\models j\imageof\lambda\in j(A)$. Hence $A\in U_j$.\smallskip

$U_j\models\xitemof{defbility:x-gen-sc-5}$:  Suppose that
$\vec{A}:=\seqof{A_\xi}{\xi\in\lambda}\in\uniV$ is \st\ $A_\xi\in U_j$, i.e.
\begin{xitemize}
\xitem[defbility:x-gen-sc-5-2-0] 
$j\imageof{\lambda}\in j(A_\xi)$
\end{xitemize}
for all $\xi<\lambda$.

By elementarity, we have
\begin{xitemize}
\xitem[defbility:x-gen-sc-5-3] $j(\triangle_{\xi\in\lambda}A_\xi)=\setof{x\in\Pkl{j(\kappa)}{j(\lambda)}^M}{
\forall\eta\in x\,(x\in j(\vec{A}(\eta)))}$
\end{xitemize}

For $\eta\in j\imageof{\lambda}$, there is $\eta_0\in\lambda$ \st\ $\eta=j(\eta_0)$. Thus 
\begin{xitemize}
\xitem[defbility:x-gen-sc-5-4] $j(\vec{A})(\eta)=j(\vec{A})(j(\eta_0))
\obecause{=}{}{by elementarity}
\ubecause{j(\vec{A}(\eta_0))}{}{$=j(A_{\eta_0})$}
\obecause{\ni}{}{\xitemof{defbility:x-gen-sc-5-2-0}}j\imageof{\lambda}.$
\end{xitemize}
By \xitemof{defbility:x-gen-sc-5-3} and \xitemof{defbility:x-gen-sc-5-4}, it follows that
$j\imageof{\lambda}\in j(\triangle_{\xi\in\lambda}A_\xi)$, and thus
$\triangle_{\xi\in\lambda}A_\xi\in U_j$.\\
\qedofClaim\qedskip

It follows that there is $\condp\in\genG$ \st
\begin{xitemize}
\xitem[def:x-gen-sc-0] $\condp\forces{\poP}{\mbox{there is a }\uniV\mbox{-normal ultrafilter on }
  \calP^\uniV(\Pkl{}{}^\uniV)}.$
\end{xitemize}
Since $\poP\restr\condp\in\calP$ by the assumption on $\calP$, we obtain the desired 
situation for $\lambda$ by replacing $\poP$ with $\poP\restr\condp$. 
\qedskip

``$\Leftarrow$'':\quad
Let $\lambda\geq\kappa$ and let $\poP$ be a $\LT\mu$-closed \po\ with 
a $(\uniV,\poP)$-generic $\genG$ and $\uniV$-normal ultrafilter $U\in\uniV[\genG]$ on
$\calP^\uniV(\Pkl{}{}^\uniV)$. 

Let
\begin{xitemize}
\xitem[] 
  $\calW:=\setof{f\in \uniV}{\mapping{f}{\Pkl{}{}^\uniV}{\uniV}}$  
\end{xitemize}
\memo{$\calP^V(\Pkl{}{}^\uniV)$ → $\Pkl{}{}^\uniV$}

\begin{xitemize}
\xitem[defbility:x-gen-sc-5-5] 
  For $f$, $g\in\calW$, $f\sim_U\,g$ $\ \ :\Leftrightarrow$\ \ 
  $\setof{x\in\Pkl{}{}^\uniV}{f(x)=g(x)}\in U$;\\[\jot]
  \phantom{For $f$, $g\in\calW$, }$f\in_U\,g$\ \  $:\Leftrightarrow$\ \ 
  $\setof{x\in\Pkl{}{}^\uniV}{f(x)\in g(x)}\in U$.
\end{xitemize}

$\sim_U$ is a congruence relation to $\in_U$. Thus may consider $\in_U$ as a binary 
relation on $\calW/{\sim}_U$  and simply write 
\begin{xitemize}
\xitem[]  
	$f/{\sim}_U\ \in_U\,g/{\sim}_U$ $\ \ :\Leftrightarrow$\ \  $f\in_U\,g$.\,\footnotemark
\end{xitemize}
\footnotetext{\phantomsection\label{scott}Here we 
  apply the common trick to handle the equivalence classes 
  by defining 
  \begin{xitemize}
  \item[] $f/{\sim}_U:=\setof{g\in\calW}{
    \begin{array}[t]{@{}l}
      g\sim_U f\mbox{ and }g\mbox{ is of minimal 
        $\in$-rank}\\ \mbox{among elements of $\calW$ with this property}}
    \end{array}
    $
  \end{xitemize}
  to make each equivalence class $f/{\sim}_U$ a set. 
} 

Let $\mapping{i_U}{\uniV}{\calW/{\sim}_U}$ be defined by
\begin{xitemize}
\xitem[defbility:x-gen-sc-6] 
  $i_U(a):=const_a/{\sim}_U$
\end{xitemize}
for $a\in\uniV$ where $const_a$ denote the function on $\Pkl{}{}^\uniV$ whose value is 
constantly $a$. \L o\'s's Theorem holds: 
\begin{Claim}
\Label{Cl-gen-sc-0-0}
For any formula $\varphi=\varphi(\variables{x}{n-1})$ in $\Lin$ (the language of \ZF), and
$\variables{f}{n-1}\in\calW$,
we have\quad
$\pairof{\calW/{\sim}_U,\in_U}\models\varphi(f_0/{\sim}_U\ctentenc f_{n-1}/{\sim}_U)$,\quad if and only 
if\\
$\setof{x\in\Pkl{}{}^\uniV}{\uniV\models\varphi(f_0(x)\ctentenc f_{n-1}(x))}\in U$. 
\end{Claim}
\prfofClaim By induction on $\varphi$. \qedofClaim\qedskip

By \Claimabove, the class mapping $i_U$ above is an elementary embedding of $\uniV$ into
$\pairof{\calW/{\sim}_U,\in_U}$. 

\begin{Claim}
  \Label{Cl-gen-sc-0-1}
  $\in_U$ is \assertof{i} an extensional, \assertof{ii} well-founded and \assertof{iii} set-like 
  relation on $\calW/{\sim}_U$. 
\end{Claim}
\prfofClaim \assertof{i}:\quad The extensionality of $\in_U$ follows from the elementarity of $i_U$.

\assertof{ii}:\quad Assume, toward a contradiction, that there is a 
sequence $\seqof{f_n}{n\in\omega}$ in $\calW$ \st\ $f_{n+1}\in_U f_n$ for all $n\in\omega$. 
By the definition of $\in_U$, this means that
$A_n:=\setof{x\in\Pkl{}{}^\uniV}{f_{n+1}(x)\in f_n(x)}\in U$ for all $n\in\omega$. 
Since $\poP$ does not add any new $\omega$-sequence, $\seqof{f_n}{n\in\omega}\in\uniV$. 
Thus, we also have $\seqof{A_n}{n\in\omega}\in\uniV$. 
By \Lemmaof{def:P-gen-sc-2},\,\assertof{1}, it follows that $\bigcap_{n\in\omega}A_n\in U$. 
For an element $x$ of this intersection, we have
\begin{xitemize}
\item[] $f_0(x)\ni f_1(x)\ni f_2(x)\ni f_3(x)\ni\,\cdots$
\end{xitemize}
by definition of $A_n$'s. 
This is a contradiction. \smallskip

\assertof{iii}:\quad  Let $f\in\calW$ be arbitrary, and let
$S:=\bigcup_{x\in\Pkl{}{}^\uniV}f(x)$. Then, by \L o\,s's Theorem, we have
\begin{xitemize}
\item[] $\setof{g/{\sim}_U}{g/{\sim}_U\,\in_U f/{\sim}_U}
\subseteq\setof{g/{\sim}_U}{\mapping{g}{\Pkl{}{}^\uniV}{S}}$
\end{xitemize}
The right side of the inclusion is clearly a set. \qedofClaim\qedskip

Let $\mapping{\mu_U}{\pairof{\calW/{\sim}_U, \in_U}}{\pairof{M,\in}}$
be the Mostowski-collapse, and let $\mapping{[\cdot]_U}{\calW}{M}$; $f\mapsto[f]_U:=\mu_U(f/{\sim}_U)$.

\L \'os's Theorem (\Claimof{Cl-gen-sc-0-0}) translates to the following: 
\begin{Claim}
  \Label{Cl-gen-sc-1} For any formula $\varphi=\varphi(\variables{x}{n-1})$ in $\Lin$ (the language of \ZF), and
  $\variables{f}{n-1}\in\calW$, we have\quad 
  $M\models\varphi([f_0]_U\ctentenc [f_{n-1}]_U)$,\quad if and only if\\
  $\setof{x\in\Pkl{}{}^\uniV}{\uniV\models\varphi(f_0(x)\ctentenc f_{n-1}(x))}\in U$. \qedofClaimmark
\end{Claim}

Let 
\begin{xitemize}
\item[] $\elembed{j_U}{\uniV}{M}$; $a\mapsto [a]_U:=\mu_U(i_U(a))=[const_a]_U$. 
\end{xitemize}
We show that 
$\elembed{j_U}{\uniV}{M}$ 
is a $\lambda$-generically supercompact embedding for $\kappa$. 

\begin{Claim}\Label{Cl-gen-sc-2}
  \wassertof{1}\quad $j_U(\xi)=\xi$ for all $\xi\in\kappa$.\smallskip

  \wassert{2} $j_U\imageof\lambda\in M$. \smallskip

  \wassert{3} $j_U(\kappa)>\lambda$.

\end{Claim}
\prfofClaim \assertof{1}: Note that $j_U(\xi)=\mu_U(i_U(\xi))=[const_\xi]_U$. Thus, for 
$\xi<\kappa$ and $f\in\calW$, 
\begin{xitemize}
\item[] $[f]_U\in j_U(\xi)\ \ \Leftrightarrow\ \ [f]_U\in [const_\xi]_U\\[\jot] 
  \ubecause{\Leftrightarrow}{}{\Claimof{Cl-gen-sc-1}}\ \setof{x\in\Pkl{}{}^\uniV}{f(x)\in 
    \ubecause{\xi}{}{\qquad\qquad$=const_\xi(x)$}}\in U\\[\jot]
  \ubecause{\Leftrightarrow}{}{\qquad by \Lemmaof{def:P-gen-sc-2},\,\assertof{2} and \xitemof{defbility:x-gen-sc-4}}
  \setof{x\in\Pkl{}{}^\uniV}{f(x)=\ubecause{\eta^*}{}{\qquad\qquad$=const_{\eta^*}(x)
      $}}\in U\mbox{ for some }\eta^*\in\xi
  \\[2\jot]
  \ubecause{\Leftrightarrow}{}{\Claimof{Cl-gen-sc-1}}\ [f]_U=j_U(\eta^*)\mbox{ for some }\eta^*\in\xi.
  $
\end{xitemize}
Thus, by induction on $\xi<\kappa$, we obtain $j_U(\xi)=\xi$ for all $\xi<\kappa$. 
\smallskip

\assertof{2}: We show that $[id_{\Pkl{}{}^\uniV}]_U=j_U\imageof{\lambda}$.

For an arbitrary $f\in\calW$ 
\begin{xitemize}
\item[] $[f]_U\in[id_{\Pkl{}{}^\uniV}]_U
  \ubecause{\Leftrightarrow}{}{by \Claimof{Cl-gen-sc-1}}
  \setof{x\in\Pkl{}{}^\uniV}{f(x)\in\ubecause{x}{}{\qquad\qquad$=id_{\Pkl{}{}^\uniV}(x)$}}\in U\\[\jot]
  \ubecause{\Leftrightarrow}{}{by \Lemmaof{def:P-gen-sc-2},\,\assertof{2}}  
  \setof{x\in\Pkl{}{}^\uniV}{f(x)=\ubecause{\xi^*}{}{\qquad\qquad$=const_{\xi^*}(x)$}}\in U
  \mbox{\/ for some }\xi^*<\lambda\\[\jot]
  \ubecause{\Leftrightarrow}{}{by \Claimof{Cl-gen-sc-1}}  
  [f]_U=j_U(\xi^*)\mbox{ for some }\xi^*<\lambda.$\ifarxived%
\fi
\end{xitemize}

\assertof{3}: 
We have 
\begin{xitemize}
\item[]  
$M\modelof{\otp([\id_{\Pkl{}{}^\uniV}]_U)<j(\kappa)}$
\end{xitemize}
by \L o\'s's Theorem 
(\Claimof{Cl-gen-sc-1}) since
$\setof{z\in\Pkl{}{}^\uniV}{\otp(x)<\ubecause{\!\kappa\!}{}{$=const_\kappa(x)$}}=\Pkl{}{}^\uniV\in U$. 

On the other hand:
\begin{xitemize}
\item[]
$M\modelof{\otp([\id_{\Pkl{}{}^\uniV}]_U)\smash{\obecause{=}{}{by \assertof{2}}}\lambda\ }$.
\qedofClaim\vspace{-1.8ex}
\end{xitemize}
\qedof{\Thmof{def:P-gen-sc-1}}
\qedskip

Note that the proof of \Claimof{Cl-gen-sc-0-1} relies on the condition on $\calP$ that no
$\poP\in\calP$ adds any new $\omega$-sequence ground model sets. Note also that the argument using the
fact that the well-foundedness of a relation is $\Delta_1$ is irrelevant here since the relation $\in_U$ is 
not in the ground model. 

Thus, the proof of \Thmof{def:P-gen-sc-1} cannot simply 
be applied to the generic supercompactness by a class of \pos\ $\calP$ whose elements might 
add new $\omega$-sequences of ground model sets. 







By \Thmof{def:P-gen-sc-1} we obtain another characterization of generic supercompactness by a 
$\calP$ as in \Thmof{def:P-gen-sc-1}:

\begin{Cor}
\Label{def:P-gen-sc-2-0} Suppose that $\calP$ is a class of \pos\ \st\ no $\poP\in\calP$ adds any 
new $\omega$-sequence of ground model sets, and $\calP$ is closed \wrt\ restriction. Then,  
the following are equivalent:\smallskip\smallskip

\wassert{a} $\kappa$ is generically supercompact by $\calP$.\smallskip

\wassert{b} For any $\lambda\geq\kappa$, there is a $\poP\in\calP$ \st\ 
\begin{xitemize}
\item[] 
  $\forces{\poP}{\mbox{there is a }\uniV\mbox{-normal ultrafilter on }
  \calP^\uniV(\Pkl{}{}^\uniV)}.$
\end{xitemize}

\wassert{c} For any $\lambda\geq\kappa$, there is a $\poP\in\calP$ \st\ for any
$(\uniV,\poP)$-generic $\genG$, there are classes $j$, $M\subseteq\uniV[\genG]$ \st\
$\elembed{j}{\uniV}{M}\subseteq\uniV[\genG]$; $\crit(j)=\kappa$; $j(\kappa)>\lambda$ and
$j\imageof{\lambda}\in M$. \qed
\end{Cor}

For a class $\calP$ of \pos\ which may contain \pos\ adding a new $\omega$ sequence of 
ground model sets, we have to modify the argument above to obtain the following theorem 
which also implies the definability of generic supercompactness by $\calP$. 

We shall call a $\uniV$-normal ultrafilter $U$ on $\calP^\uniV(\Pkl{}{}^\uniV)$ {\it steep}
if $\in_U$ defined as in \xitemof{defbility:x-gen-sc-5-5} is well-founded.

\begin{Thm}
\Label{def:P-gen-sc-3} Suppose that $\calP$ is a class of \pos\ \st\ $\calP$ is closed \wrt\ restriction. Then,  
the following are equivalent:\smallskip\smallskip

\wassert{a} $\kappa$ is generically supercompact by $\calP$.\smallskip

\wassert{b} For any regular $\lambda\geq\kappa$, there is a $\poP\in\calP$ \st\ 
\begin{xitemize}
\item[] 
  $\forces{\poP}{\mbox{there is a steep }\uniV\mbox{-normal ultrafilter on }
  \calP^\uniV(\Pkl{}{}^\uniV)}.$
\end{xitemize}

\wassert{c} For any $\lambda\geq\kappa$, there is a $\poP\in\calP$ \st\ for any
$(\uniV,\poP)$-generic $\genG$, there are classes $j$, $M\subseteq\uniV[\genG]$ \st\
$\elembed{j}{\uniV}{M}\subseteq\uniV[\genG]$, $\crit(j)=\kappa$, $j(\kappa)>\lambda$, 
and
\quad$j\imageof{\lambda}\in M$. 
\end{Thm}
\noindent
\prfof{\Thmof{def:P-gen-sc-3}} A slight modification the proof of \Thmof{def:P-gen-sc-1} 
will do: it is enough 
to show that, for $U_j$ in the proof of ``$\Rightarrow$'' of \Thmof{def:P-gen-sc-1}, the 
relation $\in_{U_j}$ defined in \xitemof{defbility:x-gen-sc-5-5} 
is well-founded. This follows from the next Claim:
\begin{Claim}
\Label{Cl-gen-sc-3} In $\uniV[\genG]$, the class mapping
\begin{xitemize}
\xitem[] $\mapping{\iota}{\calW/{\sim}_{U_j}}{\uniV[\genG]}$;\quad $f/{\sim}_{U_j}\ \ \mapsto\ \  j(f(j\imageof\lambda))$
\end{xitemize}
is well-defined, and it is an embedding of $\pairof{\calW/{\sim}_{U_j},\,\in_{U_j}}$ into $\pairof{\uniV[\genG],\in}$. 
\end{Claim}
\prfofClaim For $f$, $g\in\calW$, we have
\begin{xitemize}
\item[]  $f/{\sim}_{U_j}\ \sim_{U_j}\ g/{\sim_{U_j}}\ \ 
\ubecause{\Leftrightarrow}{}{by the definition \xitemof{defbility:x-gen-sc-5-5} of $\sim_{U_j}$}\ \ 
\setof{x\in\Pkl{}{}^\uniV}{f(x)=g(x)}\in U_j
\ubecause{\Leftrightarrow}{}{by the definition \xitemof{defbility:x-gen-sc-5-2-a} of $U_j$ }\ \ \\[2\jot]
j(\setof{x\in\Pkl{}{}^\uniV}{f(x)=g(x)})\ni j\imageof{\lambda}
\ \ \Leftrightarrow\ \ \ubecause{j(f)(j\imageof\lambda)}{}{\qquad\qquad$=\iota(f/{\sim}_{U_j})$}
=\ubecause{j(g)(j\imageof\lambda)}{}{\qquad\qquad$=\iota(g/{\sim}_{U_j})$}$.
\end{xitemize}
This shows the well-definedness and the injectivity of $\iota$. 

Similarly we can show 
\begin{xitemize}
\item[]  
  $f/{\sim}_{U_j}\ \in_{U_j}\ g/{\sim_{U_j}}\ \ \Leftrightarrow\ \ 
  \ubecause{j(f)(j\imageof{\lambda})}{}{\qquad\qquad$=\iota(f/{\sim}_{U_j})$}
  \in\ubecause{j(g)(\imageof{\lambda})}{}{\qquad\qquad$=\iota(g/{\sim}_{U_j})$}$. \ifarxived\vspace{-3.6ex}\else\vspace{-2zh}\fi
\end{xitemize}
\qedofClaim\\
\qedof{\Thmof{def:P-gen-sc-3}}

\section{Sufficiently large initial segment of elementary embeddings}
\Label{defbility:initial}
In this section, we prove a characterization of Laver-generic supercompactness from which 
the first-order definability of this notion follows.  

\begin{Lemma}
\Label{def:P-is-0-0}
Suppose that $\poP$ is a \po\ (in $\uniV$), and $\genG$ 
a $(\uniV,\poP)$-generic set. Suppose that $j$, $M\subseteq\uniV[\genG]$ are \st\ 
$\elembed{j}{\uniV}{M\subseteq\uniV[\genG]}$.

Then, for a cardinal $\theta$ (in $\uniV$), 
have:\quad $\elembed{j\restr\calH(\theta)^\uniV}{\calH(\theta)^\uniV}{\calH(j(\theta))^M}$.
\end{Lemma}
\prf For any $\Lin$-formula $\varphi=\varphi(\variables{x}{k-1})$ and
$\variables{u}{k-1}\in\calH(\theta)^\uniV$, we have
\begin{xitemize}
\item[]
$\calH(\theta)^\uniV\models\varphi(\variables{u}{k-1})\ \ \Leftrightarrow\ \ 
\uniV\modelof{\calH(\theta)^\uniV\models\varphi(\variables{u}{k-1})}\\[\jot]
\mbox{}\ifarxived
\hspace{-0.72ex}
\else\hspace{-0.4zw}\fi\ubecause{\Leftrightarrow}{}{\qquad\qquad\qquad by elementarity of $j$}\ \ 
M\modelof{\calH(j(\theta))^M\models\varphi(j(u_0)\ctentenc j(u_{k-1}))}\\[3\jot]
\Leftrightarrow\ \ \calH(j(\theta))^M\models\varphi(j(u_0)\ctentenc j(u_{k-1}))$. \qedofLemma
\end{xitemize}

Note that, in the Lemma above, $\calH(j(\theta))^M$ is transitive since $M$ is transitive. 

\begin{Lemma}
\Label{def:P-is-0-1} Suppose that $\poP$ is a \po\ (in $\uniV$), and $\genG$ 
a $(\uniV,\poP)$-generic set. Suppose further that $\theta$ is a cardinal in $\uniV$ and 
$j_0, N\in\uniV[\genG]$ be \st\ $N$ is transitive 
and $\elembed{j_0}{\calH(\theta)^\uniV}{N}$.

Let $N_0=\bigcup j_0\imageof{\calH(\theta)^\uniV}$. Then, we have:\smallskip

\wassert{1} $N_0$ is transitive.\smallskip

\wassert{2} \assertof{i}\/ $N_0\prec N$, \quad\assertof{ii}\/ $j_0\imageof{\calH(\theta)}\subseteq N_0$, and
\quad \assertof{iii}\/ $\elembed{j_0}{\calH(\theta)^\uniV}{N_0}$.\smallskip

\wassert{3} For any $b\in N_0$, there is $a\in\calH(\theta)^\uniV$ \st\ $b\in j_0(a)$. \smallskip

\wassert{4} If $\theta_0<\theta$ is \st\ $\calH(\theta_0)^\uniV\in\calH(\theta)^\uniV$ then
$\calH(j_0(\theta_0))^N\subseteq N_0$. 
\end{Lemma}
\prf \assertof{1}: Suppose that $b\in N_0$ and $c\in b$. We have to show that $c\in N_0$. 

Let $a\in\calH(\theta)^\uniV$ be \st\ $b\in j_0(a)$. Let $a^*=\trcl(a)$. Then
$a^*\in\calH(\theta)^{\uniV}$. Since
$\calH(\theta)^\uniV\models a^*\mbox{ is transitive and }a\subseteq a^*$, we have 
\begin{xitemize}
\item[] $M\models j_0(a^*)\mbox{ is transitive and }j(a)\subseteq j(a^*)$
\end{xitemize}
by elementarity. Since $N$ is transitive, $j_0(a^*)$ is really transitive. Since
$c\in b\in j_0(a^*)$, it follows that
$c\in j_0(a^*)\subseteq\bigcup j_0\imageof\calH(\theta)^\uniV=N_0$. 
\smallskip

\assertof{2}, \assertof{i}: We check that $N_0$ satisfies Vaught's criterion.

Suppose that $b_1\ctentenc b_n\in N_0$ and $\varphi(\variables{x}{n})$ is an $\Lin$-formula 
\st\ 
\begin{xitemize}
\xitem[def:x-is-2-a-0]  
$N\models\exists x\varphi(x,b_1\ctentenc b_n)$.
\end{xitemize}

We have to show that there is 
$b\in N_0$ \st\ $N\models\varphi(b,b_1\ctentenc b_n)$.

Let $a_i\in\calH(\theta)^\uniV$ for $i\in n+1\setminus 1$ be \st\ $b_i\in j_0(a_i)$ for all
$i\in n+1\setminus 1$. 
Then we have
\begin{xitemize}
\xitem[def:x-is-2-a] $\calH(\theta)^\uniV\models\exists x\forall y_1\in a_1\cdots\forall y_n\in a_n
\,\Big(\exists y\varphi(y,y_1\ctentenc,y_n)$\\
\phantom{$\calH(\theta)^\uniV\models\exists x\forall y_1\in a_1\cdots\forall y_n\in a_n
\,\Big($}
$\rightarrow\,\exists y\in x\,\varphi(y,y_1\ctentenc y_n)\Big).
$
\end{xitemize}

Let $a\in\calH(\theta)^\uniV$ be a witness of \xitemof{def:x-is-2-a}. That is,
\begin{xitemize}
\item[]  $\calH(\theta)^\uniV\models\forall y_1\in a_1\cdots\forall y_n\in a_n
\,\Big(\exists y\varphi(y,y_1\ctentenc,y_n)$\\
\phantom{$\calH(\theta)^\uniV\models\forall y_1\in a_1\cdots\forall y_n\in a_n
\,\Big($}
$
\rightarrow\,\exists y\in a\,\varphi(y,y_1\ctentenc y_n)\Big).$
\end{xitemize}
By elementarity, it follows that 
\begin{xitemize}
\xitem[def:x-is-2-a-1]  $N\models\forall y_1\in j_0(a_1)\cdots\forall y_n\in j_0(a_n)
\,\Big(\exists y\varphi(y,y_1\ctentenc,y_n)$\\
\phantom{$N\models\forall y_1\in j_0(a_1)\cdots\forall y_n\in j_0(a_n)
\,\Big($}
$
\rightarrow\,\exists y\in j_0(a)\,\varphi(y,y_1\ctentenc y_n)\Big).$
\end{xitemize}

By \xitemof{def:x-is-2-a-1} and \xitemof{def:x-is-2-a-0}, there is
$b\in j_0(a)\subseteq\bigcup j_0\imageof{\calH(\theta)^\uniV}=N_0$ \st\ 
\begin{xitemize}
\item[]  
$N\models\varphi(b,b_1\ctentenc b_n)$.
\end{xitemize}

\assertof{2},\,\assertof{ii}: Suppose that $a\in\calH(\theta)^\uniV$. Then
$\ssetof{a}\in\calH(\theta)^\uniV$ and
$j_0(a)\in\ssetof{j_0(a)}=j_0(\ssetof{a})\subseteq\bigcup j_0\imageof{\calH(\theta)^\uniV}=N_0$.
\smallskip

\assertof{2},\,\assertof{iii}: This follows from 
\assertof{2},\,\assertof{i},\,\assertof{ii}. 
\smallskip

\assertof{3}: This is clear by definition of $N_0$. \smallskip

\assertof{4}: Suppose that $\theta_0<\theta$ is \st\
$\calH(\theta_0)^\uniV\in\calH(\theta)^\uniV$. Let $a=\calH(\theta_0)^\uniV$. By 
elementarity, $N\models j_0(a)\mbox{ is }\calH(j(\theta_0))$. Thus
$j_0(a)=\calH(j(\theta_0))^N$ and $j_0(a)\in N_0$ by \assertof{2}, \assertof{ii}. By 
\assertof{1}, it follows that $\calH(j(\theta_0))^N\subseteq N_0$.
\qedofLemma

\begin{Prop}
\Label{def:P-is-1}
Suppose that $\poP$ is a \po\ (in $\uniV$) and $\genG$ a $(\uniV,\poP)$-generic filter. 
Suppose further that $\theta$ is a regular cardinal and
$\elembed{j_0}{\calH(\theta)^\uniV}{N}$ for a transitive set $N\in\uniV[\genG]$ \st, 
\begin{xitemize}
\xitem[def:x-is-3] 
$\poP\in\calH(\theta)^\uniV$; and, \smallskip
\xitem[def:x-is-3-0] 
 for any $b\in N$, there is $a\in\calH(\theta)^\uniV$ \st\
$b\in j_0(a)$. 
\end{xitemize}
Then there are $j$, $M\subseteq\uniV[\genG]$ \st\ 
\begin{xitemize}
\xitem[def:x-is-4]  $\elembed{j}{\uniV}{M\subseteq\uniV[\genG]}$,
\xitem[def:x-is-5] $N\subseteq M$ and $j\restr\calH(\theta)^\uniV=j_0$.
\end{xitemize}
\end{Prop}
\prf
We mainly work in $\uniV[\genG]$. Let
\begin{xitemize}
\xitem[def:x-is-6]
$\calF:=\setof{f\in\uniV}{\mapping{f}{\dom(f)}{\uniV},\dom(f)\in\calH(\theta)^\uniV}$, and
\xitem[def:x-is-7] $\Pi:=\setof{\pairof{f,a}}{f\in\calF,\,a\in j_0(\dom(f))}$. 
\end{xitemize}

For $\pairof{f,a}$, $\pairof{g,b}\in\Pi$, let
\begin{xitemize}
\xitem[def:x-is-8]  $\pairof{f,a}\sim\pairof{g,b}$\ \ $:\Leftrightarrow$\ \ $\pairof{a,b}\in j_0(S_{f(x)=g(y)})$,\\
where
$S_{f(x)=g(y)}:=\setof{\pairof{x,y}}{x\in\dom(f),\,y\in\dom(g),\,f(x)=g(x)}$; 
and
\xitem[def:x-is-9]  $\pairof{f,a}\mathrel{E}\pairof{g,b}$\ \ $:\Leftrightarrow$\ \
$\pairof{a,b}\in j_0(S_{f(x)\varin g(y)})$,\\
where
$S_{f(x)\varin g(y)}:=\setof{\pairof{x,y}}{x\in\dom(f),\,y\in\dom(g),\,f(x)\in g(x)}$. 
\end{xitemize}

\begin{Claim}\Label{def:Cl-0}\wassertof{1}\quad $\sim$ is an equivalence relation on\/ $\Pi$.
\smallskip

\wassert{2} $\sim$ is a congruence relation to $E$.
\end{Claim}
\prfofClaim
\assertof{1}: Clearly $\sim$ is reflective and symmetric. We show that $\sim$ is 
transitive. Suppose that $\pairof{f,a}$, $\pairof{g,b}$, $\pairof{h,c}\in\Pi$, 
$\pairof{f,a}\sim\pairof{g,b}$ and $\pairof{g,b}\sim\pairof{h,c}$. By the definition
\xitemof{def:x-is-8}, we have $\pairof{a,b}\in j_0(S_{f(x)=g(y)})$ and
$\pairof{b,c}\in j_0(S_{g(y)=h(z)})$. Thus
\begin{xitemize}
\item[]  $\pairof{a,c}\in j_0(S_{f(x)=g(y)})\circ j_0(S_{g(y)=h(z)})
\ubecause{=}{}{by elementarity of $j_0$} j_0(S_{f(x)=g(y)}\circ S_{g(y)=h(z)})\\
\mbox{}\!\!\!\!\ubecause{\subseteq}{}{\qquad\qquad\qquad\qquad\qquad\qquad\qquad\quad
by $S_{f(x)=g(y)}\circ S_{g(y)=h(z)}\subseteq S_{f(x)=h(z)}$ and elementarity}
j_0(S_{f(y)=h(z)})$.
\end{xitemize}
This shows that $\pairof{f,a}\sim\pairof{h,c}$.\smallskip

\assertof{2}: Suppose $\pairof{f_0,a_0}$, $\pairof{f_1,a_1}$, $\pairof{g,b}\in\Pi$,\quad
$\pairof{f_0,a_0}\sim \pairof{f_1,a_1}$, \quad and \\
$\pairof{f_0,a_0}\mathrel{E}\pairof{g,b}$. 
Then

\begin{xitemize}
\item[]  $\pairof{a_1,b}\in j_0(S_{f_1(x_1)=f_0(x_0)})\circ j_0(S_{f_0(x_0)\varin g(y)})=
j_0(S_{f_1(x_1)=f_0(x_0)}\circ S_{f_0(x_0)\varin g(y)})\\[\jot]\subseteq j_0(S_{f_1(x_1)\varin g(y)})$.
\end{xitemize}
Thus $\pairof{f_1,a_1}\mathrel{E}\pairof{g,b}$.

Similarly, we can show that, for $\pairof{f,a}$, $\pairof{g_0,b_0}$,
$\pairof{g_1,b_1}\in\Pi$, $\pairof{g_0,b_0}\sim \pairof{g_1,b_1}$ and
$\pairof{f,a}\mathrel{E}\pairof{g_0,b_0}$ implies
$\pairof{f,a}\mathrel{E}\pairof{g_1,b_1}$. Since $\sim$ is a equivalence relation by 
\assertof{1}, it follows that $\sim$ is a congruence relation to $E$. 
\qedofClaim
\qedskip

Let $\Pi/{\sim}$ be the class of the equivalence classes (in the sense of \footnoteof{scott}) of
$\sim$.  We denote the equivalence class of $\pairof{f,a}\in\Pi$ modulo $\sim$
by $\pairof{f,a}/{\sim}$. 
For simplicity, we denote the binary relation on $\Pi/{\sim}$ corresponding to $E$ 
also by $E$. Thus, $\pairof{f,a}/{\sim}\mathrel{E}\pairof{g,b}/{\sim}$ $:\Leftrightarrow$
$\pairof{f,a}\mathrel{E}\pairof{g,b}$. 

Generalizing the notation we already used in \xitemof{def:x-is-8} and \xitemof{def:x-is-9}, we let
\begin{xitemize}
\item[\qquad\quad\ ] $S_{\varphi(f_0(x_0)\ctentenc f_{n-1}(x_n-1))}\\
\qquad:=\setof{\pairof{\variables{u}{n-1}}\in\uniV}{
\begin{array}[t]{@{}l}
u_0\in\dom(f_0)\ctentenc u_{n-1}\in\dom(f_{n-1}),\\
\uniV\models\varphi(f_0(u_0)\ctentenc f_{n-1}(u_{n-1}))\quad}
\end{array}
$
\end{xitemize}
for each $\Lin$-formula $\varphi=\varphi(\variables{x}{n-1})$. 

We have the following ``\,\L o\'s's Theorem\,'' for $\pairof{\Pi/{\sim}, E}$.

\begin{Claim}
\Label{def:Cl-is-a} For any $\Lin$-formula $\varphi=\varphi(\variables{x}{n-1})$ and
$\pairof{f_0,a_0}\ctentenc\pairof{f_{n-1},a_{n-1}}\in\Pi$, we have
\begin{xitemize}
\item[\ \ ]
$\pairof{\Pi/{\sim}, E}\models\varphi(\pairof{f_0,a_0}/{\sim}\ctentenc\pairof{f_{n-1},a_{n-1}}/{\sim})$
\\[\jot]
$\Leftrightarrow\ \ \pairof{\variables{a}{n-1}}\in j_0(S_{\varphi(f_0(x_0)\ctentenc f_{n-1}(x_{n-1}))}).$
\end{xitemize}
\end{Claim}
\prfofClaim
By induction on $\varphi$. If $\varphi$ is atomic, the claim follows from the definitions \xitemof{def:x-is-8}
and \xitemof{def:x-is-9} of $\sim$ and $E$. 

The induction step for ``$\varphi=\neg\varphi_0$'' is trivial.

Suppose $\varphi=\varphi(\variables{x}{n-1})$, 
$\varphi=\varphi_0\lor\varphi_1$, and 
$\pairof{f_0,a_0}\ctentenc\pairof{f_{n-1},a_{n-1}}\in\Pi$. Note that
\begin{xitemize}
\xitem[def:x-is-10]
$S_{\varphi(f_0(x_0)\ctentenc f_{n-1}(x_{n-1}))}=S_{\varphi_0(f_0(x_0)\ctentenc f_{n-1}(x_{n-1}))}\cup 
S_{\varphi_1(f_0(x_0)\ctentenc f_{n-1}(x_{n-1}))}$.
\end{xitemize}

We have
\begin{xitemize}
\item[] $\pairof{\Pi/{\sim}, E}\models\varphi(\pairof{f_0,a_0}/{\sim}\ctentenc\pairof{f_{n-1},a_{n-1}}/{\sim})$\\[\jot]
\ $\Leftrightarrow$\ \ $\pairof{\Pi/{\sim}, E}\models\varphi_0(\pairof{f_0,a_0}/{\sim}\ctentenc\pairof{f_{n-1},a_{n-1}}/{\sim})$
\\ \phantom{\ $\Leftrightarrow$\ \ }
or $\pairof{\Pi/{\sim}, E}\models\varphi_1(\pairof{f_0,a_0}/{\sim}\ctentenc\pairof{f_{n-1},a_{n-1}}/{\sim})$\\[2\jot]
\ifarxived
$\obecause{\Leftrightarrow}{1.26ex}{\qquad by induction hypothesis}$\ \ 
\else
$\obecause{\Leftrightarrow}{0.7zh}{\qquad by induction hypothesis}$\ \ 
\fi
$\pairof{\variables{a}{n-1}}\in j_0(S_{\varphi_0(f_0(x_0)\ctentenc f_{n-1}(x_{n-1}))})$
\\ \phantom{\ $\Leftrightarrow$\ \ }
or $\pairof{\variables{a}{n-1}}\in j_0(S_{\varphi_1(f_0(x_0)\ctentenc f_{n-1}(x_{n-1}))})$\\[\jot]
\ $\Leftrightarrow$\ \ 
$\pairof{\variables{a}{n-1}}\in j_0(S_{\varphi_0(f_0(x_0)\ctentenc f_{n-1}(x_{n-1}))})\cup 
j_0(S_{\varphi_1(f_0(x_0)\ctentenc f_{n-1}(x_{n-1}))})$\\[2\jot]
\ifarxived
$\obecause{\Leftrightarrow}{1.26ex}{\qquad\qquad by elementarity of $j$ and \xitemof{def:x-is-10} }$\ \ 
\else
$\obecause{\Leftrightarrow}{0.7zh}{\qquad\qquad by elementarity of $j$ and \xitemof{def:x-is-10} }$\ \ 
\fi
$\pairof{\variables{a}{n-1}}\in j_0(S_{\varphi(f_0(x_0)\ctentenc f_{n-1}(x_{n-1}))})$. 
\end{xitemize}

Finally, suppose $\varphi=\exists x\varphi_0(x,x_1\ctentenc x_{n-1})$ and
$\pairof{f_1,a_1}\ctentenc\pairof{f_{n-1},a_{n-1}}\in\Pi$.  

If
$\pairof{\Pi/{\sim}, E}\models\varphi(\pairof{f_1,a_1}/{\sim}\ctentenc\pairof{f_{n-1},a_{n-1}}/{\sim})$, 
then there is $\pairof{f,a}\in\Pi$ \st\
$\pairof{\Pi/{\sim}, E}\models\varphi_0(\pairof{a,f}/{\sim}, \pairof{f_1,a_1}/{\sim}\ctentenc\pairof{f_{n-1},a_{n-1}}/{\sim})$. 
By induction hypothesis, it follows that
$\pairof{a,a_1\ctentenc a_{n-1}}\in j_0(S_{\varphi_0(f(x_0),f_1(x_1)\ctenten)})$. Thus, by 
elementarity and by the definition of $S_{\varphi(\cdots)}$,
$\pairof{a_1\ctentenc a_{n-1}}\in j_0(S_{\varphi(f_1(x_1)\ctentenc f_{n-1}(x_{n-1}))})$. 

Conversely, assume that 
$\pairof{a_1\ctentenc a_{n-1}}\in j_0(S_{\varphi(f_1(x_1)\ctentenc f_{n-1}(x_{n-1}))})$.  
Let $d=\dom(f_1)\times\cdots\times\dom(f_{n-1})$. Note that $d\in\calH(\theta)^\uniV$.

Let
$f\in\uniV$ with $\mapping{f}{d}{\uniV}$ be defined by 
\begin{xitemize}
\item[]  $f(\pairof{\variables{u}{n-1}})=\left\{\ 
\begin{array}{@{}l}
\mbox{some }u\in\uniV\mbox{ \st\ }\calH(\theta)^\uniV\models\varphi_0(u,\variables{u}{n-1}),\\
\phantom{\emptyset,}\qquad\qquad\qquad\qquad\mbox{ if there is such }u\in\uniV\,;\\[2\jot]
\emptyset,\qquad\qquad\qquad\qquad\mbox{ otherwise.}
\end{array}
\right.
$
\end{xitemize}

We have
\begin{xitemize}
\item[]  $\calH(\theta)^\uniV\models{}
\begin{array}[t]{@{}l}
  \forall x_1\cdots\forall x_{n-1}\Big(\pairof{x_1\ctentenc x_{n-1}}\in S_{\varphi(f_1(x_1)\ctenten)}\\[\jot]
  \qquad\rightarrow\exists x\,(\pairof{x,x_1\ctentenc x_{n-1}}\in S_{\varphi_0(f(x),f_1(x_1)\ctenten)}
  \Big).
\end{array}
$
\end{xitemize}

By elementarity, it follows that 
\begin{xitemize}
\item[]  $N\models{}
\begin{array}[t]{@{}l}
  \forall x_1\cdots\forall x_{n-1}\Big(\pairof{x_1\ctentenc x_{n-1}}\in j_0(S_{\varphi(f_1(x_1)\ctenten)})\\[\jot]
  \qquad\rightarrow\exists x\,(\pairof{x,x_1\ctentenc x_{n-1}}\in j_0(S_{\varphi_0(f(x),f_1(x_1)\ctenten)})
  \Big).
\end{array}
$
\end{xitemize}

Hence, there is $a\in N$ \st\
$\pairof{a,a_1\ctentenc a_{n-1}}\in j_0(S_{\varphi_0(f(x),f_1(x_1)\ctenten)})$. By induction 
hypothesis, it follows that
\begin{xitemize}
\item[]
$\pairof{\Pi/{\sim}, E}\models\varphi_0(\pairof{a,f}/{\sim}, \pairof{f_1,a_1}/{\sim}\ctentenc\pairof{f_{n-1},a_{n-1}}/{\sim})$. 
\end{xitemize}
Thus
$\pairof{\Pi/{\sim}, E}\models\varphi(\pairof{f_1,a_1}/{\sim}\ctentenc\pairof{f_{n-1},a_{n-1}}/{\sim})$. 
\qedofClaim
\qedskip

For $u\in\uniV$, let $\mapping{f_u}{1}{\uniV}$ be defined by $f_u(\emptyset)=u$. 
Let $\mapping{ i }{\uniV}{\Pi/{\sim}}$ be defined by $ i (u)=\pairof{f_u,\emptyset}/{\sim}$. 
\begin{Claim}
\Label{def:Cl-is-0} $ i $ is an elementary embedding of $\pairof{\uniV,\in}$ into
$\pairof{\Pi/{\sim}, E}$.
\end{Claim}
\prfofClaim
Suppose that $\varphi=\varphi(\variables{x}{n-1})$ is an $\Lin$-formula and
$\variables{u}{n-1}\in\uniV$. Then we have
\begin{xitemize}
\item[] $\pairof{\Pi/{\sim},E}\models\varphi( i (u_0)\ctentenc i (u_{n-1}))$\\[\jot]
$\ubecause{\Leftrightarrow}{}{by \Claimof{def:Cl-is-a}}\ \ \pairof{\emptyset,\emptyset\ctentenc\emptyset}\in
\begin{array}[t]{@{}l}
  j_0(S_{\varphi(f_{u_0}(x_0)\ctentenc f_{u_{n-1}}(x_{n-1}))})\\[2\jot]
  \ifarxived
  \mbox{}\hspace{-0.7em}\obecause{=}{1.44ex}{\qquad\qquad\qquad by definition of $S_{\varphi(\cdots)}$}
  \else
  \mbox{}\hspace{-0.7zw}\obecause{=}{0.8zh}{\qquad\qquad\qquad by definition of $S_{\varphi(\cdots)}$}
  \fi
  j_0(\setof{\pairof{\variables{x}{n-1}}}{\uniV\models\varphi(f_{u_0}(x_0)\ctentenc f_{u_{n-1}}(x_{n-1}))})\\[\jot]
  =\left\{\ 
  \begin{array}{@{}ll}
  \emptyset, &\mbox{if }\uniV\not\models\varphi(u_0\ctentenc u_{n-1});\\[\jot]
  \ssetof{\pairof{\emptyset\ctentenc\emptyset}}, &\mbox{if }\uniV\models\varphi(u_0\ctentenc u_{n-1}).
\end{array}
\right.\\[\jot]
\end{array}\\[\jot]
\ \Leftrightarrow\ \ \uniV\models\varphi(u_0\ctentenc u_{n-1}).
$
\end{xitemize}\ifarxived\vspace{-1.8ex}\else\vspace{-1zh}\fi

\qedofClaim
\begin{Claim}
\Label{def:Cl-is-1} \wassertof{1}\quad $E$ is well-founded.\smallskip

\wassert{2} $E$ is set like.
\end{Claim}
\prfofClaim \assertof{1}: Suppose not and let $\pairof{f_n,b_n}\in\Pi$, $n\in\omega$ (in $V[\genG]$) be \st\ 
\begin{xitemize}
\xitem[def:x-is-11]  $\pairof{f_0,b_0}\mathrel{\reflectbox{$E$}}\pairof{f_1,b_1}\mathrel{\reflectbox{$E$}}
                     \pairof{f_2,b_2}\mathrel{\reflectbox{$E$}}\ \cdots$.
\end{xitemize}

Let $\utf_n$, $n\in\omega$ be $\poP$-names of $f_n$, $n\in\omega$ (note that we can choose $\utf_n$,
$n\in\omega$ \st\ 
$\seqof{\utf_n}{n\in\omega}\in\uniV$), and let
\begin{xitemize}
\xitem[def:x-is-12]
$\calQ:=\setof{\pairof{\condp,n,u}}{{}
\begin{array}[t]{@{}l}
\condp\in\poP,\,n\in\omega,\,u\in\calH(\theta)^\uniV,\\[\jot]
\condp\mbox{ decides }\utf_n,
\mbox{ and }\condp\forces{\poP}{u\varin\dom(\utf_n)}\ }.
\end{array}
$
\end{xitemize}
By \xitemof{def:x-is-3} and since $\theta$ is regular, we have $\calQ\in\calH(\theta)^\uniV$. 

For $\pairof{\condp_0,n_0,u_0}$, $\pairof{\condp_1,n_1,u_1}\in\calQ$, let 
\begin{xitemize}
\item[] $\pairof{\condp_0,n_0,u_0}\sqsubset\pairof{\condp_1,n_1,u_1}$\ \ $:\Leftrightarrow$\ \ 
        $\condp_0\leq_\poP\condp_1$,\quad$n_0=n_1+1$,\\
        \phantom{$\pairof{\condp_0,n_0,u_0}\sqsubset\pairof{\condp_1,n_1,u_1}$\ \ $:\Leftrightarrow$\ \ }%
        and\quad $\condp_0\forces{\poP}{\utf_{n_0}(u_0)\varin\utf_{n_1}(u_1)}$. 
\end{xitemize}

In $\uniV[\genG]$, let $\seqof{\condp_n}{n\in\omega}$ be a descending sequence in $\genG$ 
\wrt\ $\leq_\poP$ 
\st\ each $\condp_n$ decides $\utf_n$ to be $f_n$. 
\begin{Subclaim}
$\seqof{\pairof{j_0(\condp_n),n,b_n}}{n\in\omega}$ is a descending sequence in
$j_0(\pairof{\calQ,\sqsubset})$ \wrt\ $j_0(\sqsubset)$.
\end{Subclaim}
\prfofClaim For $n\in\omega$, we have to show that 
\begin{xitemize}
\item[]
$\pairof{j_0(\condp_{n+1}),n+1,b_{n+1}}\mathrel{j_0(\sqsubset)} \pairof{j_0(\condp_0), n,b_n}$
\end{xitemize}
holds. By the choice of $\condp_n$'s, we have $\condp_{n+1}\leq_\poP\condp_n$,
$\condp_{n+1}\forces{\poP}{\utf_{n+1}=f_{n+1}}$, and
$\condp_{n}\forces{\poP}{\utf_{n}=f_{n}}$. Thus we have
\begin{xitemize}
\xitem[def:x-is-12-a]
$\condp_{n+1}\forces{\poP}{\utf_{n+1}=f_{n+1}\ \land\ \utf_{n}=f_{n}}$. 
\end{xitemize}

It follows that 
\begin{xitemize}
\item[]  $\sqsubset {}
\begin{array}[t]{@{}l}
  \supseteq\setof{\pairof{\pairof{\condp_{n+1},n+1,u},\pairof{\condp_{n},n,v}}}{
  \condp_{n+1}\forces{\poP}{\utf_{n+1}(u)\varin\utf_n(v)}}\\[\jot]
  \!\!\ubecause{=}{}{by \xitemof{def:x-is-12-a}}
  \setof{\pairof{\pairof{\condp_{n+1},n+1,u},\pairof{\condp_{n},n,v}}}{
  f_{n+1}(u)\in f_n(v) }\\[\jot]
  \!\!\ubecause{=}{}{\qquad\qquad by the definition of $S_{\threedotsh\varin\threedotsh}$ in \xitemof{def:x-is-9}}
  \setof{\pairof{\pairof{\condp_{n+1},n+1,u},\pairof{\condp_{n},n,v}}}{
  \pairof{u,v}\in S_{f_{n+1}(x_0)\varin f_n(x_1)}}.
\end{array}
$
\end{xitemize}

Thus 

\begin{xitemize}
\item[]  
$j_0(\sqsubset)\ {}
\begin{array}[t]{@{}l}
  \supseteq\setof{(\pairof{\pairof{j_0(\condp_{n+1}),n+1,u},\pairof{j_0(\condp_{n}),n,v}})}{
  \pairof{u,v}\in j_0(S_{f_{n+1}(x_0)\varin f_n(x_1)})}\\[\jot]
  \ni\pairof{\pairof{j_0(\condp_{n+1}),n+1,b_{n+1}}, \pairof{j_0(\condp_0), n,b_n}}. 
  \end{array}$\vspace{-4ex}
  \end{xitemize}
\qedofSubclaim
\qedskip

Since being well-founded is $\Delta_1$, it follows that
$N\modelof{j_0(\pairof{\calQ,\sqsubset})\mbox{ is not well-founded}}$. By elementarity,
it follows that 
$\calH(\theta)^\uniV\modelof{\pairof{\calQ,\sqsubset}\mbox{ is not well-founded}}$. 
However, if 
$\seqof{\pairof{\condq_n, k_n, u_n}}{n\in\omega}$ is a descending sequence in
$\pairof{\calQ,\sqsubset}$, then we would have
\begin{xitemize}
\item[] $g_{k_0}(u_0)\ni g_{k_1}(u_1)\ni g_{k_2}(u_2)\ni\ \cdots$\ 
\end{xitemize}
where $g_{k_n}$, for each $n\in\omega$, is the element of $\calF$ which is decided to be $\utf_{k_n}$ \ifarxived
\vspace{-1.08ex}
\else
\vspace{-0.4zh}
\fi
by $\condp_n$. 
This is a contradiction.\smallskip

\assertof{2}: Suppose that $\pairof{f,a}$, $\pairof{g,b}\in\Pi$ and 
\begin{xitemize}
\xitem[] $\pairof{f,a}\mathrel{E}\pairof{g,b}$.  
\end{xitemize}
Let $\mapping{f_0}{\dom(f)}{\bigcup g\imageof{\dom(g)}\cup\ssetof{\infty}}$, where $\infty$ 
is a set \st\ $\infty\not\in g\imageof{\dom(g)}$, be defined by
\begin{xitemize}
\item[]  $f_0(u)=\left\{\,
\begin{array}{@{}ll}
  f(u), &\mbox{if }f(u)\in\bigcup g\imageof{\dom(g)};\\[\jot]
  \infty, &\mbox{otherwise}
\end{array}
\right.$
\end{xitemize}
for all $u\in\dom(f)$. By the definition of $f_0$, we have
$S_{f(x_0)\varin g(x_1)}=S_{f_0(x_0)\varin g(x_1)}$. Thus we have 
\begin{xitemize}
\xitem[]  $\pairof{f,a}\sim\pairof{f_0,a}$.
\end{xitemize}

This implies that

\begin{xitemize}
\item[] $\setof{\pi\in\Phi/{\sim}}{\pi\mathrel{E}\pairof{g,b/{\sim}}}$\\[\jot]
$\subseteq\setof{\pairof{f,a}/{\sim}}{
\begin{array}[t]{@{}l}
\dom(f)\in\calH(\theta)^\uniV, \\[\jot]
\mapping{f}{\dom(f)}{\bigcup g\imageof{\dom(g)}\cup\ssetof{\infty}},\ a\in j_0(\dom(f))}
\end{array}
$.
\end{xitemize}
The right side of the inclusion is clearly a set. 
\qedofClaim
\qedskip

$\pairof{\Pi/{\sim},E}$ is extensional by \Claimof{def:Cl-is-0}. Hence, by 
\Claimof{def:Cl-is-1}, there is the Mostowski collapse
\begin{xitemize}
\item[]  $\mapping{m}{\pairof{\Pi/{\sim},E}}{\pairof{\uniV[\genG],\in}}$. 
\end{xitemize}

Let $M:=m\imageof{\Pi/{\sim}}$ and $j:=m\circ i$. By \Claimof{def:Cl-is-0}, we have
\begin{xitemize}
\item[]  $\elembed{j}{\uniV}{M\subseteq\uniV[\genG]}$. 
\end{xitemize}

Note that, for $a\in\calH(\theta)^\uniV$, 
\begin{xitemize}
\xitem[def:x-is-12-0]  $j(a)=m\circ i(a)=m(\pairof{f_a,\emptyset}/{\sim})$. 
\end{xitemize}

For each $b\in N$, let $d_b\in\calH(\theta)^\uniV$ be \st\ $b\in j_0(d_b)$. We can 
always find such $d_b$ by \xitemof{def:x-is-3-0}. Let
\begin{xitemize}
\item[]  $\mapping{\iota}{N}{\Pi/{\sim}}$;\quad $b\mapsto\pairof{\id_{d_b},b}/{\sim}$. 
\end{xitemize}
\begin{Claim}
\Label{def:Cl-is-2} $\iota$ is an embedding of $\pairof{N,\in}$ into $\pairof{\Pi/{\sim},E}$, and 
$\iota\imageof{N}$ is a full initial segment of $\Pi/{\sim}$ \wrt\ $E$. In particular, for 
any $b\in N$, we have $m(\iota(b))=m(\pairof{id_{d_b},b}/{\sim})=b$.
\end{Claim}
\prfofClaim
Note that 
\begin{xitemize}
\xitem[def:x-is-13]  
$j_0(\id_{d_b})=\id_{j_0(d_b)}$ 
\end{xitemize}
by elementarity.

For $b$, $c\in N$\ifarxived\vspace{-1.8ex}\else\vspace{-1zh}\fi
\begin{xitemize}
\item[]  
$\iota(b)\mathrel{E}\iota(c)\ \ \ubecause{\Leftrightarrow}{}{by definition of $\iota$}\ \ 
\pairof{id_{d_b},b}\mathrel{E}\pairof{id_{d_c},c}\ \ 
\obecause{\Leftrightarrow}{}{by the definition \xitemof{def:x-is-9} of $E$}
\ubecause{j_0(\id_{d_b})(b)}{}{\qquad\qquad$=b$, by \xitemof{def:x-is-13}}\in 
\obecause{j_0(\id_{d_c})(c)}{}{\qquad\qquad$=c$, by \xitemof{def:x-is-13}}.
$
\end{xitemize}

Suppose that $\pairof{f,a}/{\sim}\mathrel{E}\pairof{\id_{d_b},b}=\iota(b)$ 
for $\pairof{f,a}\in\Pi$. This means that 
\begin{xitemize}
\item[] $j_0(f)(a)\in j_0(\id_{d_b})(b)\ubecause{=}{}{by \xitemof{def:x-is-13}}b$.
\end{xitemize}
Let $c:=j_0(f)(a)$. Then we have $c\in b\in N$. Since $N$ is transitive it follows that
$c\in N$. By the definition \xitemof{def:x-is-8} of $\sim$, we have
\begin{xitemize}
\item[]  
$\iota(c)=\pairof{id_{d_c},c}/{\sim}=\pairof{f,a}/\sim$. \qedofClaim
\end{xitemize}

Together with the previous Claim, the following Claim shows that our $j$ and $M$ are as desired:
\begin{Claim}
\Label{def:Cl-is-3} $j\restr\calH(\theta)^\uniV=j_0$. 
\end{Claim}
\prfofClaim
Suppose that $a\in\calH(\theta)^\uniV$. We show that $j(a)=j_0(a)$.

Note that 
$j(a)=m(\pairof{f_a,\emptyset}/{\sim})$. 
For $b:=j_0(a)$, we have $\pairof{f_a,\emptyset}\sim\pairof{id_{d_b},b}$
by \xitemof{def:x-is-8}. It follows that $j(a)=m(\pairof{\id_{d_b},b}/{\sim})
\ubecause{=}{}{by \Claimof{def:Cl-is-2}}b=j_0(a). $
\qedofClaim\\[-2.4ex]
\qedofProp

\begin{Thm}
\Label{def:P-is-2} Suppose that $\calP$ is an iterable class of \pos. Then \tfae:
\smallskip

\wassert{a} $\kappa$ is L-g supercompact for $\calP$.
\smallskip

\wassert{b} For any $\lambda$, and for any $\poP\in\calP$, there is a $\poP$-name $\utpoQ$ 
with $\forces{\poP}{\utpoQ\varin\poP}$ \st\
\begin{xitemize}
\item[]
$\forcessanf{\poP\ast\utpoQ}{\parbox[t]{0.9\textwidth}{
there are a regular cardinal $\theta$, a transitive set $N$, and a mapping $j_0$ \\
\st\\[2\jot]
\assertof{1}\quad $\elembed{j_0}{\calH(\theta)^\uniV}{N}$,\\[\jot]
\assertof{2}\quad $\crit(j)=\kappa$, $\theta,\,j(\kappa)>\lambda$,\\[\jot]
\assertof{3}\quad for any $b\varin N$, there is $a\varin\calH(\theta)^\uniV$ \st\ $b\varin j_0(a)$\,\\[\jot]
\assertof{4}\quad $\poP\ast\utpoQ$, $\utilde{\genH}\in N$, and\\
\assertof{5}\quad $j\imageof\lambda\in N$\ {\rm''}.
}}
$
\end{xitemize}
\end{Thm}
\prf ``\assertof{a} $\Rightarrow$ \assertof{b}'': By \Lemmaof{def:P-is-0-0} and \Lemmaof{def:P-is-0-1}. 

``\assertof{b} $\Rightarrow$ \assertof{a}'': By \Propof{def:P-is-1}.\qedofThm

\phantomsection

\end{document}